\input amssym.def
\input psfig
\input epsf

\let \blskip = \baselineskip
\parskip=1.2ex plus .2ex minus .1ex
\let\s=\span
\tabskip 20pt
\tolerance = 1000
\pretolerance = 50
\newcount\itemnum
\itemnum = 0
\overfullrule = 0pt

\def\title#1{\bigskip\centerline{\bigbigbf#1}}
\def\author#1{\bigskip\centerline{\bf #1}\smallskip}
\def\address#1{\centerline{\it#1}}
\def\abstract#1{\vskip1truecm{\narrower\noindent{\bf Abstract.} #1\bigskip}}

\def\sp{\bigskip}
\def\nosp{\vskip -\the\blskip plus 1pt minus 1pt}

\def\br{\hfil\break} 
\def\ti{\br \hglue \the \parindent}

\def\ce#1{\LP\centerline{#1}}

\def\skipit#1{}
\def\mdag{\raise 3pt\hbox{\dag}}

\def\IP{\par\hang}
\def\XP{\par\noindent\hang}
\def\LP{\par\noindent}
\def\BP[#1]{\par\item{[#1]}}
\def\SH#1{\sp\vskip\parskip\leftline{\bigbf #1}\nobreak}

\def\TH#1{\sp\XP{\bf THEOREM\ \shead#1}}
\def\LM#1{\sp\XP{\bf LEMMA\ \shead#1}}

\def\CO#1{\sp\XP{\bf COROLLARY\ \shead#1}}

\def\PF{\LP{\bf Proof:\ }}
\def\NX{\advance\itemnum by 1 \sp\LP {\bf \shead \the\itemnum.\ }}
\def\qed{\null\nobreak\hfill\hbox{${\vrule width 5pt height 6pt}$}\par\sp}

\def\cart{\>\hbox{${\vcenter{\vbox{
    \hrule height 0.4pt\hbox{\vrule width 0.4pt height 4.5pt
    \kern4pt\vrule width 0.4pt}\hrule height 0.4pt}}}$}\>}
\def\bxmu{\>\hbox{${\vcenter{\vbox {
    \hrule height 0.4pt\hbox{\vrule width 0.4pt height 4pt
    \hskip -1.3pt\lower 1.8pt\hbox{$\times$}\negthinspace\vrule width 0.4pt}
    \hrule height 0.4pt}}}$}\>}

\def\lin#1{\hbox to #1true in{\hrulefill}}


\def\foot#1{\raise 6pt \hbox{#1} \kern -3pt}

\def\fig #1 #2 #3 #4 #5 {\sp \ce{ {\epsfbox[#1 #2 #3 #4]{figs/#5.ps}} }}

\def\gpic#1{#1 \sp\ce{\box\graph} \medskip} 
\def\ccol{\hfil##\hfil&}

\def\rcol{\hfil##&}

\def\TS#1{\medskip \ce{\vbox{\offinterlineskip\tabskip=1em
                  \halign{\strut#1}}}\medskip}


\def\DM{{\it Discrete Math.{}}}

\def\TAMS{{\it Trans.\ Amer.\ Math.\ Soc.{}}}

\def\AMM{{\it Amer.\ Math.\ Monthly}}

\def\EJC{{\it Europ.\ J.\ Comb.{}}}

			\def\dlt{\delta}
\def\eps{\epsilon}    
	 	
  \def\SG{\Sigma}






\def\({\left(}	\def\){\right)}


\def\CH#1#2{{{#1}\choose{#2}}}

\def\FR#1#2{{#1 \over #2}}

\def\FL#1{\left\lfloor{#1}\right\rfloor}

\def\SE#1#2#3{\sum_{#1 = #2} ^ {#3}}

\def\PE#1#2#3{\prod_{#1 = #2} ^ {#3}}

\def\VEC#1#2#3{#1_{#2},\ldots,#1_{#3}}

\def\MAP#1#2#3{#1\colon\;#2\to#3}
\def\SET#1:#2{\{#1\colon\;#2\}}


		
\def\C#1{\left | #1 \right |}    
	\def\un#1{\underline{#1}}

\def\AND{\ {\rm and}\ }

\def\mod{\,{\rm mod}\,}





\magnification=\magstep1
\vsize=9.0 true in
\hsize=6.5 true in
\headline={\hfil\ifnum\pageno=1\else\folio\fi\hfil}
\footline={\hfil\ifnum\pageno=1\folio\else\fi\hfil}

\parindent=20pt
\baselineskip=12pt
\parskip=.5ex  

\def\shead{ }

\font\bigbf = cmb10 scaled \magstep1

\font\bigmit = cmmi10 scaled \magstep1
\font\bigbigbf = cmb10 scaled \magstep2


\title{THE BRICKLAYER PROBLEM}
\title{AND THE STRONG CYCLE LEMMA}
\author{Hunter S. Snevily}
\address{University of Idaho, Moscow, ID 83844-1103, {snevily@uidaho.edu}}
\author{Douglas B. West} 
\address{University of Illinois, Urbana, IL 61801-2975, {west@math.uiuc.edu}}
\vfootnote{}{\br
   Running head: {BRICKLAYER PROBLEM AND CYCLE LEMMA} \br 
   Keywords: Cycle Lemma, Catalan numbers\br
   Written July 1996, revised May 1997.
}


\def\Ckf{\FR1{k+1}\CH{2k}k}

\def\Cqn{\FR1{qn+1}\CH{(q+1)n}n}

\SH
{1. THE PROBLEM}
One of the most familiar brand names in the world of toys is $Lego^{TM}$;
the name immediately brings interlocking building blocks to mind.
Interlocking blocks fit together in restricted ways.  Typically, a block of
order $q$ (length $q+1$) has $q+1$ protrusions on its top and $q+1$
indentations on its
bottom, so that the indentations on the bottom of one block can lock on to
the protrusions on the top of another.  Here $q$ is a positive integer, and
the width and height of a block are unimportant.  For $q=1$ and $q=2$,
Figure 1 shows the $q$ ways in which one block can sit on top of two others.
\gpic{
\expandafter\ifx\csname graph\endcsname\relax \csname newbox\endcsname\graph\fi
\expandafter\ifx\csname graphtemp\endcsname\relax \csname newdimen\endcsname\graphtemp\fi
\setbox\graph=\vtop{\vskip 0pt\hbox{%
    \special{pn 8}%
    \special{pa 1429 600}%
    \special{pa 1500 600}%
    \special{pa 1500 571}%
    \special{pa 1643 571}%
    \special{pa 1643 600}%
    \special{pa 1714 600}%
    \special{fp}%
    \special{pa 1429 314}%
    \special{pa 1500 314}%
    \special{pa 1500 286}%
    \special{pa 1643 286}%
    \special{pa 1643 314}%
    \special{pa 1714 314}%
    \special{fp}%
    \special{pa 1714 600}%
    \special{pa 1786 600}%
    \special{pa 1786 571}%
    \special{pa 1929 571}%
    \special{pa 1929 600}%
    \special{pa 2000 600}%
    \special{fp}%
    \special{pa 1714 314}%
    \special{pa 1786 314}%
    \special{pa 1786 286}%
    \special{pa 1929 286}%
    \special{pa 1929 314}%
    \special{pa 2000 314}%
    \special{fp}%
    \special{pa 1429 600}%
    \special{pa 1429 314}%
    \special{fp}%
    \special{pa 2000 600}%
    \special{pa 2000 314}%
    \special{fp}%
    \special{pa 2000 600}%
    \special{pa 2071 600}%
    \special{pa 2071 571}%
    \special{pa 2214 571}%
    \special{pa 2214 600}%
    \special{pa 2286 600}%
    \special{fp}%
    \special{pa 2000 314}%
    \special{pa 2071 314}%
    \special{pa 2071 286}%
    \special{pa 2214 286}%
    \special{pa 2214 314}%
    \special{pa 2286 314}%
    \special{fp}%
    \special{pa 2286 600}%
    \special{pa 2357 600}%
    \special{pa 2357 571}%
    \special{pa 2500 571}%
    \special{pa 2500 600}%
    \special{pa 2571 600}%
    \special{fp}%
    \special{pa 2286 314}%
    \special{pa 2357 314}%
    \special{pa 2357 286}%
    \special{pa 2500 286}%
    \special{pa 2500 314}%
    \special{pa 2571 314}%
    \special{fp}%
    \special{pa 2000 600}%
    \special{pa 2000 314}%
    \special{fp}%
    \special{pa 2571 600}%
    \special{pa 2571 314}%
    \special{fp}%
    \special{pa 1714 314}%
    \special{pa 1786 314}%
    \special{pa 1786 286}%
    \special{pa 1929 286}%
    \special{pa 1929 314}%
    \special{pa 2000 314}%
    \special{fp}%
    \special{pa 1714 29}%
    \special{pa 1786 29}%
    \special{pa 1786 0}%
    \special{pa 1929 0}%
    \special{pa 1929 29}%
    \special{pa 2000 29}%
    \special{fp}%
    \special{pa 2000 314}%
    \special{pa 2071 314}%
    \special{pa 2071 286}%
    \special{pa 2214 286}%
    \special{pa 2214 314}%
    \special{pa 2286 314}%
    \special{fp}%
    \special{pa 2000 29}%
    \special{pa 2071 29}%
    \special{pa 2071 0}%
    \special{pa 2214 0}%
    \special{pa 2214 29}%
    \special{pa 2286 29}%
    \special{fp}%
    \special{pa 1714 314}%
    \special{pa 1714 29}%
    \special{fp}%
    \special{pa 2286 314}%
    \special{pa 2286 29}%
    \special{fp}%
    \special{pa 0 1457}%
    \special{pa 71 1457}%
    \special{pa 71 1429}%
    \special{pa 214 1429}%
    \special{pa 214 1457}%
    \special{pa 286 1457}%
    \special{fp}%
    \special{pa 0 1171}%
    \special{pa 71 1171}%
    \special{pa 71 1143}%
    \special{pa 214 1143}%
    \special{pa 214 1171}%
    \special{pa 286 1171}%
    \special{fp}%
    \special{pa 286 1457}%
    \special{pa 357 1457}%
    \special{pa 357 1429}%
    \special{pa 500 1429}%
    \special{pa 500 1457}%
    \special{pa 571 1457}%
    \special{fp}%
    \special{pa 286 1171}%
    \special{pa 357 1171}%
    \special{pa 357 1143}%
    \special{pa 500 1143}%
    \special{pa 500 1171}%
    \special{pa 571 1171}%
    \special{fp}%
    \special{pa 571 1457}%
    \special{pa 643 1457}%
    \special{pa 643 1429}%
    \special{pa 786 1429}%
    \special{pa 786 1457}%
    \special{pa 857 1457}%
    \special{fp}%
    \special{pa 571 1171}%
    \special{pa 643 1171}%
    \special{pa 643 1143}%
    \special{pa 786 1143}%
    \special{pa 786 1171}%
    \special{pa 857 1171}%
    \special{fp}%
    \special{pa 0 1457}%
    \special{pa 0 1171}%
    \special{fp}%
    \special{pa 857 1457}%
    \special{pa 857 1171}%
    \special{fp}%
    \special{pa 857 1457}%
    \special{pa 929 1457}%
    \special{pa 929 1429}%
    \special{pa 1071 1429}%
    \special{pa 1071 1457}%
    \special{pa 1143 1457}%
    \special{fp}%
    \special{pa 857 1171}%
    \special{pa 929 1171}%
    \special{pa 929 1143}%
    \special{pa 1071 1143}%
    \special{pa 1071 1171}%
    \special{pa 1143 1171}%
    \special{fp}%
    \special{pa 1143 1457}%
    \special{pa 1214 1457}%
    \special{pa 1214 1429}%
    \special{pa 1357 1429}%
    \special{pa 1357 1457}%
    \special{pa 1429 1457}%
    \special{fp}%
    \special{pa 1143 1171}%
    \special{pa 1214 1171}%
    \special{pa 1214 1143}%
    \special{pa 1357 1143}%
    \special{pa 1357 1171}%
    \special{pa 1429 1171}%
    \special{fp}%
    \special{pa 1429 1457}%
    \special{pa 1500 1457}%
    \special{pa 1500 1429}%
    \special{pa 1643 1429}%
    \special{pa 1643 1457}%
    \special{pa 1714 1457}%
    \special{fp}%
    \special{pa 1429 1171}%
    \special{pa 1500 1171}%
    \special{pa 1500 1143}%
    \special{pa 1643 1143}%
    \special{pa 1643 1171}%
    \special{pa 1714 1171}%
    \special{fp}%
    \special{pa 857 1457}%
    \special{pa 857 1171}%
    \special{fp}%
    \special{pa 1714 1457}%
    \special{pa 1714 1171}%
    \special{fp}%
    \special{pa 286 1171}%
    \special{pa 357 1171}%
    \special{pa 357 1143}%
    \special{pa 500 1143}%
    \special{pa 500 1171}%
    \special{pa 571 1171}%
    \special{fp}%
    \special{pa 286 886}%
    \special{pa 357 886}%
    \special{pa 357 857}%
    \special{pa 500 857}%
    \special{pa 500 886}%
    \special{pa 571 886}%
    \special{fp}%
    \special{pa 571 1171}%
    \special{pa 643 1171}%
    \special{pa 643 1143}%
    \special{pa 786 1143}%
    \special{pa 786 1171}%
    \special{pa 857 1171}%
    \special{fp}%
    \special{pa 571 886}%
    \special{pa 643 886}%
    \special{pa 643 857}%
    \special{pa 786 857}%
    \special{pa 786 886}%
    \special{pa 857 886}%
    \special{fp}%
    \special{pa 857 1171}%
    \special{pa 929 1171}%
    \special{pa 929 1143}%
    \special{pa 1071 1143}%
    \special{pa 1071 1171}%
    \special{pa 1143 1171}%
    \special{fp}%
    \special{pa 857 886}%
    \special{pa 929 886}%
    \special{pa 929 857}%
    \special{pa 1071 857}%
    \special{pa 1071 886}%
    \special{pa 1143 886}%
    \special{fp}%
    \special{pa 286 1171}%
    \special{pa 286 886}%
    \special{fp}%
    \special{pa 1143 1171}%
    \special{pa 1143 886}%
    \special{fp}%
    \special{pa 2286 1457}%
    \special{pa 2357 1457}%
    \special{pa 2357 1429}%
    \special{pa 2500 1429}%
    \special{pa 2500 1457}%
    \special{pa 2571 1457}%
    \special{fp}%
    \special{pa 2286 1171}%
    \special{pa 2357 1171}%
    \special{pa 2357 1143}%
    \special{pa 2500 1143}%
    \special{pa 2500 1171}%
    \special{pa 2571 1171}%
    \special{fp}%
    \special{pa 2571 1457}%
    \special{pa 2643 1457}%
    \special{pa 2643 1429}%
    \special{pa 2786 1429}%
    \special{pa 2786 1457}%
    \special{pa 2857 1457}%
    \special{fp}%
    \special{pa 2571 1171}%
    \special{pa 2643 1171}%
    \special{pa 2643 1143}%
    \special{pa 2786 1143}%
    \special{pa 2786 1171}%
    \special{pa 2857 1171}%
    \special{fp}%
    \special{pa 2857 1457}%
    \special{pa 2929 1457}%
    \special{pa 2929 1429}%
    \special{pa 3071 1429}%
    \special{pa 3071 1457}%
    \special{pa 3143 1457}%
    \special{fp}%
    \special{pa 2857 1171}%
    \special{pa 2929 1171}%
    \special{pa 2929 1143}%
    \special{pa 3071 1143}%
    \special{pa 3071 1171}%
    \special{pa 3143 1171}%
    \special{fp}%
    \special{pa 2286 1457}%
    \special{pa 2286 1171}%
    \special{fp}%
    \special{pa 3143 1457}%
    \special{pa 3143 1171}%
    \special{fp}%
    \special{pa 3143 1457}%
    \special{pa 3214 1457}%
    \special{pa 3214 1429}%
    \special{pa 3357 1429}%
    \special{pa 3357 1457}%
    \special{pa 3429 1457}%
    \special{fp}%
    \special{pa 3143 1171}%
    \special{pa 3214 1171}%
    \special{pa 3214 1143}%
    \special{pa 3357 1143}%
    \special{pa 3357 1171}%
    \special{pa 3429 1171}%
    \special{fp}%
    \special{pa 3429 1457}%
    \special{pa 3500 1457}%
    \special{pa 3500 1429}%
    \special{pa 3643 1429}%
    \special{pa 3643 1457}%
    \special{pa 3714 1457}%
    \special{fp}%
    \special{pa 3429 1171}%
    \special{pa 3500 1171}%
    \special{pa 3500 1143}%
    \special{pa 3643 1143}%
    \special{pa 3643 1171}%
    \special{pa 3714 1171}%
    \special{fp}%
    \special{pa 3714 1457}%
    \special{pa 3786 1457}%
    \special{pa 3786 1429}%
    \special{pa 3929 1429}%
    \special{pa 3929 1457}%
    \special{pa 4000 1457}%
    \special{fp}%
    \special{pa 3714 1171}%
    \special{pa 3786 1171}%
    \special{pa 3786 1143}%
    \special{pa 3929 1143}%
    \special{pa 3929 1171}%
    \special{pa 4000 1171}%
    \special{fp}%
    \special{pa 3143 1457}%
    \special{pa 3143 1171}%
    \special{fp}%
    \special{pa 4000 1457}%
    \special{pa 4000 1171}%
    \special{fp}%
    \special{pa 2857 1171}%
    \special{pa 2929 1171}%
    \special{pa 2929 1143}%
    \special{pa 3071 1143}%
    \special{pa 3071 1171}%
    \special{pa 3143 1171}%
    \special{fp}%
    \special{pa 2857 886}%
    \special{pa 2929 886}%
    \special{pa 2929 857}%
    \special{pa 3071 857}%
    \special{pa 3071 886}%
    \special{pa 3143 886}%
    \special{fp}%
    \special{pa 3143 1171}%
    \special{pa 3214 1171}%
    \special{pa 3214 1143}%
    \special{pa 3357 1143}%
    \special{pa 3357 1171}%
    \special{pa 3429 1171}%
    \special{fp}%
    \special{pa 3143 886}%
    \special{pa 3214 886}%
    \special{pa 3214 857}%
    \special{pa 3357 857}%
    \special{pa 3357 886}%
    \special{pa 3429 886}%
    \special{fp}%
    \special{pa 3429 1171}%
    \special{pa 3500 1171}%
    \special{pa 3500 1143}%
    \special{pa 3643 1143}%
    \special{pa 3643 1171}%
    \special{pa 3714 1171}%
    \special{fp}%
    \special{pa 3429 886}%
    \special{pa 3500 886}%
    \special{pa 3500 857}%
    \special{pa 3643 857}%
    \special{pa 3643 886}%
    \special{pa 3714 886}%
    \special{fp}%
    \special{pa 2857 1171}%
    \special{pa 2857 886}%
    \special{fp}%
    \special{pa 3714 1171}%
    \special{pa 3714 886}%
    \special{fp}%
    \hbox{\vrule depth1.457in width0pt height 0pt}%
    \kern 4.000in
  }%
}%
}
\sp
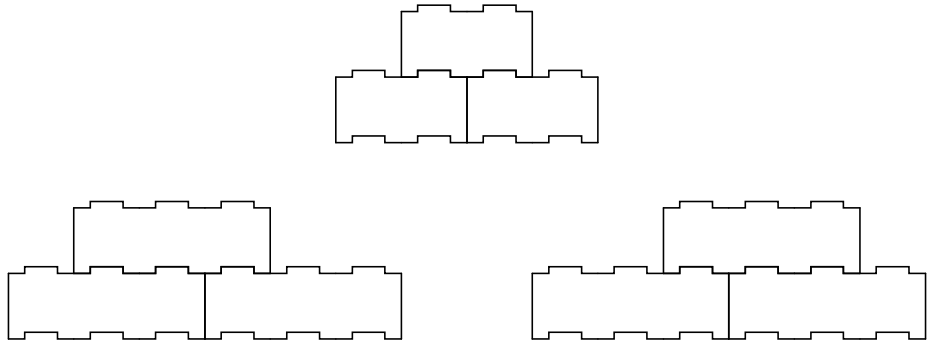
\ce{Figure 1.  Bricks of length $q+1$, for $q=1$ and $q=2$.}

In this paper, we solve a counting problem about building stacks of such
blocks, which we call {\it bricks}.  We have a linear base of length $m$
on which we can place bricks.  The bricks sitting directly on the base can
start at any integer position, as long as they fit on the base and don't
overlap.  When two bricks are contiguous, sharing a common end, we can place
another brick on top of them.  Since we want the blocks to interlock, it must
cover part of each brick below it.  The protrusions and indentations then
restrict it to $q$ possible positions, covering a positive integer amount
of each brick it rests on.  Again bricks cannot overlap.  We call a
configuration built by these rules a $q$-{\it stack}.

\IP
{\bf The bricklayer problem}.  How many different $q$-stacks can be built\br
on a base of length $m$?

For example, when $m=4$ and $q=1$, there are five $q$-stacks, including three
with one brick, one with two bricks, and one with three bricks that appears
in Figure 1.  When $m=5$ and $q=1$, there are nine $q$-stacks.  When $m=6$ and
$q=2$, there are seven $q$-stacks, two of which appear in Figure 1.

The bricklayer problem is related to a host of other problems in combinatorial
enumeration.  We solve it by establishing a bijection from the set of $q$-stacks
on a base of length $m$ to a special set of sequences of 0's and 1's.
This leads us to related counting problems because these special sequences
generalize the solutions to the famous {\it Ballot Problem}.

Consider an election that ends in a tie, with $k$ votes for each of two
candidates.  If the votes are counted in a random order, with all sequences
equally likely, we may wonder what the probability is that the first candidate
never trails.  Using 0's to designate votes for the first candidate and 1's
to designate votes for the second, we are asking for the fraction of sequences
with $k$ 1's and $k$ 0's such that every prefix (initial segment) has at least
as many 0's as 1's.  These sequences are the {\it ballot sequences} of length
$2k$.  To compute the desired probability, we count the ballot sequences and
divide by $\CH{2k}k$.  We will see shortly that the number of ballot sequences
of length $2k$ is the $k$th {\it Catalan number} $C_k$, defined by
$$C_k=\Ckf.$$

More generally, a sequence is $q$-{\it satisfying} if in each prefix the number
of 0's is at least $q$ times the number of 1's.  We will establish a bijection
between the $q$-stacks on a base of length $m$ and the $q$-satisfying sequences
of length $m$.  Thus we reduce the bricklayer problem to the counting of
$q$-satisfying sequences.

The {\it Cycle Lemma} of Dvoretzky and Motzkin [5] provides one of the many
proofs that the Catalan numbers count the ballot sequences.  With equal ease it
enables us to count the $q$-satisfying sequences of length $m$.  After doing so,
we develop the bijection to solve the Bricklayer Problem.  Subsequently, we
explore generalizations of the Cycle Lemma and applications of these 
generalizations.

\SH
{2. THE CYCLE LEMMA AND GENERALIZED CATALAN NUMBERS}
To facilitate our discussion, we introduce terminology to describe various
arrangements of 0's and 1's.  A $q$-{\it dominating} sequence (as defined in
[3]) is a sequence of 1's and 0's such that in each prefix the number of 0's is
more than $q$ times the number of 1's.  A $(k,l)$-{\it sequence} is a sequence
of $k$ 1's and $l$ 0's.  A $(k,l)$-{\it arrangement} is a cyclic arrangement of
$k$ 1's and $l$ 0's.  By this we mean that rotating the arrangement does not
change it, but we maintain a fixed direction; a flip or reversal produces
a different arrangement.

A special case of the Cycle Lemma states that every $(k,k+1)$-arrangement can be
cut in exactly one position to obtain a 1-dominating $(k,k+1)$-sequence.  The
first element of this sequence (the element after the cut) must be a 0.
Deleting this 0 yields a ballot sequence, and the process is reversible.
Thus the ballot sequences and the $(k,k+1)$-arrangements are equinumerous.
Since $k$ and $k+1$ are relatively prime, the number of $(k,k+1)$-arrangements
(and ballot sequences of length $2k$) is exactly the Catalan number
$$\FR1{2k+1}\CH{2k+1}{k+1} = \FR1{k+1}\CH{2k}k=C_k.$$

The Cycle Lemma thus ``explains'' one 0 in a $(k,k+1)$-arrangement.  Kierstead
and Trotter [9] generalized this to give combinatorial meaning to each 0.  We
present their result in Section 4 and extend it slightly.  Our task at present
is to prove the Cycle Lemma and use it to count the $q$-satisfying sequences of
length $m$.  Figure 2 illustrates the Cycle Lemma (and its proof) when
$(k,q,p)=(2,2,3)$; the underscored 0's are those that begin $2$-dominating
sequences when the arrangement is read clockwise.

\TH{1}
(The Cycle Lemma - Dvoretzky and Motzkin [5]) For $k,q,p\ge0$, every
$(k,qk+p)$-arrangement breaks to form a $q$-dominating sequence in exactly
$p$ places.
\PF
The statement is trivial for $k=0$; we procede by induction.  For $k>0$,
let $a$ be a $(k,qk+p)$-arrangement.  By the pigeonhole principle, between
some pair of the 1's in $a$ there are more than $q$ 0's.  Let $S$ be a set
of $q+1$ consecutive positions consisting of $q$ 0's followed immediately
by a 1 (illustrated by the outside arc in Figure 2).

None of the $q+1$ positions of $S$ can start a $q$-dominating sequence.
A position outside $S$ starts a $q$-dominating sequence if and only if it
starts one in the $(k-1,q(k-1)+p)$-arrangement $a'$ obtained from $a$ by
deleting $S$.  The number of $q$-dominating starting places in $a$ thus
equals the number of $q$-dominating starting places in $a'$, which by the
induction hypothesis is $p$.  \qed

\gpic{
\expandafter\ifx\csname graph\endcsname\relax \csname newbox\endcsname\graph\fi
\expandafter\ifx\csname graphtemp\endcsname\relax \csname newdimen\endcsname\graphtemp\fi
\setbox\graph=\vtop{\vskip 0pt\hbox{%
    \graphtemp=.5ex\advance\graphtemp by 0.208in
    \rlap{\kern 0.624in\lower\graphtemp\hbox to 0pt{\hss 1\hss}}%
    \graphtemp=.5ex\advance\graphtemp by 0.305in
    \rlap{\kern 0.891in\lower\graphtemp\hbox to 0pt{\hss 0\hss}}%
    \graphtemp=.5ex\advance\graphtemp by 0.552in
    \rlap{\kern 1.034in\lower\graphtemp\hbox to 0pt{\hss 0\hss}}%
    \graphtemp=.5ex\advance\graphtemp by 0.832in
    \rlap{\kern 0.984in\lower\graphtemp\hbox to 0pt{\hss 1\hss}}%
    \graphtemp=.5ex\advance\graphtemp by 1.015in
    \rlap{\kern 0.766in\lower\graphtemp\hbox to 0pt{\hss $\un0$\hss}}%
    \graphtemp=.5ex\advance\graphtemp by 1.015in
    \rlap{\kern 0.482in\lower\graphtemp\hbox to 0pt{\hss $\un0$\hss}}%
    \graphtemp=.5ex\advance\graphtemp by 0.832in
    \rlap{\kern 0.264in\lower\graphtemp\hbox to 0pt{\hss $\un0$\hss}}%
    \graphtemp=.5ex\advance\graphtemp by 0.552in
    \rlap{\kern 0.214in\lower\graphtemp\hbox to 0pt{\hss 0\hss}}%
    \graphtemp=.5ex\advance\graphtemp by 0.305in
    \rlap{\kern 0.356in\lower\graphtemp\hbox to 0pt{\hss 0\hss}}%
    \special{pn 8}%
    \special{ar 622 580 499 499 3.054345 4.799637}%
    \graphtemp=.5ex\advance\graphtemp by 0.166in
    \rlap{\kern 0.166in\lower\graphtemp\hbox to 0pt{\hss $S$\hss}}%
    \hbox{\vrule depth1.181in width0pt height 0pt}%
    \kern 1.200in
  }%
}%
}
\sp
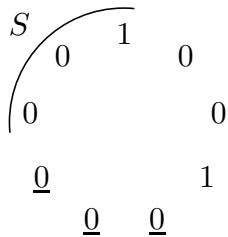
\ce{Figure 2.  2-dominating sequences from a (2,7)-arrangement.}
\skipit{\sp
\ce{1}
\ce{0\qquad\quad0}
\ce{$0\qquad\qquad\;0$}
\ce{$\un0\qquad\quad\;1$}
\ce{$\un0\quad\un0$}
}

\CO{2}
The number of $q$-satisfying $(k,qk+p-1)$-sequences and the number
of $q$-satisfying sequences of length $m$, respectively, are
$$\FR p{qk+p} \CH{(q+1)k+p-1}k \qquad{\rm and}\qquad
\SE k0{\FL{m/(q+1)}} \FR{m-(q+1)k+1}{m-k+1} \CH mk.$$
\PF
A sequence is $q$-satisfying if and only if the sequence obtained by adding a 0
at the front is $q$-dominating.  Thus the $q$-satisfying $(k,qk+p-1)$-sequences
and the $q$-dominating $(k,qk+p)$-sequences are equinumerous.  By the Cycle
Lemma, each $(k,qk+p)$-arrangement has $p$ positions that yield $q$-dominating
sequences.  Thus in each class of $(k,qk+p)$-sequences equivalent under cyclic
rotation, the fraction that are $q$-dominating is $p/(k+qk+p)$.  Since the
fraction is the same over all classes, we need not worry about periodicity.
We conclude that the number of $q$-dominating $(k,qk+p)$-sequences is
$\CH{(q+1)k+p}k p/(k+qk+p)$, which equals the formula claimed.

Each $q$-satisfying sequence of length $m$ has $k$ 1's and $m-k$ 0's, for some
$k$.  By setting $m-k=qk+p-1$, we obtain $p=m-(q+1)k+1$ and can
use the preceding formula to count these sequences.  Thus the term for $k$
in the summation is precisely the number of $q$-satisfying sequences of length
$m$ that have $k$ 1's.  \qed

Dershowitz and Zaks [3] gave two applications of the Cycle Lemma and provided
references to several proofs of it.  Many applications (and many proofs) have
been given for the special case $p=1$.  These lead to the 
{\it generalized Catalan numbers}, defined by
$$C_n^q=\Cqn.$$
By Corollary 2, $C_n^q$ is the number of $q$-satisfying $(n,qn)$-sequences.
There are other ways to generalize the Catalan numbers by introducing additional
parameters, but this is the generalization appropriate for our discussion.

For example, the generalized Catalan numbers arise in counting the rooted plane
trees with $qn+1$ leaves in which every non-leaf vertex has exactly $q+1$
children.  By ``plane tree'', we mean that the left-to-right order of children
matters.  Figure 3 shows the three such trees when $q=2$ and $n=2$.
\gpic{
\expandafter\ifx\csname graph\endcsname\relax \csname newbox\endcsname\graph\fi
\expandafter\ifx\csname graphtemp\endcsname\relax \csname newdimen\endcsname\graphtemp\fi
\setbox\graph=\vtop{\vskip 0pt\hbox{%
    \graphtemp=.5ex\advance\graphtemp by 0.563in
    \rlap{\kern 0.100in\lower\graphtemp\hbox to 0pt{\hss $\bullet$\hss}}%
    \graphtemp=.5ex\advance\graphtemp by 0.563in
    \rlap{\kern 0.350in\lower\graphtemp\hbox to 0pt{\hss $\bullet$\hss}}%
    \graphtemp=.5ex\advance\graphtemp by 0.563in
    \rlap{\kern 0.600in\lower\graphtemp\hbox to 0pt{\hss $\bullet$\hss}}%
    \graphtemp=.5ex\advance\graphtemp by 0.563in
    \rlap{\kern 0.850in\lower\graphtemp\hbox to 0pt{\hss $\bullet$\hss}}%
    \graphtemp=.5ex\advance\graphtemp by 0.563in
    \rlap{\kern 1.100in\lower\graphtemp\hbox to 0pt{\hss $\bullet$\hss}}%
    \graphtemp=.5ex\advance\graphtemp by 0.063in
    \rlap{\kern 0.850in\lower\graphtemp\hbox to 0pt{\hss $\bullet$\hss}}%
    \graphtemp=.5ex\advance\graphtemp by 0.313in
    \rlap{\kern 0.350in\lower\graphtemp\hbox to 0pt{\hss $\bullet$\hss}}%
    \special{pn 8}%
    \special{pa 100 563}%
    \special{pa 350 313}%
    \special{pa 350 563}%
    \special{pa 350 313}%
    \special{pa 600 563}%
    \special{fp}%
    \special{pa 350 313}%
    \special{pa 850 63}%
    \special{pa 850 563}%
    \special{pa 850 63}%
    \special{pa 1100 563}%
    \special{fp}%
    \graphtemp=.5ex\advance\graphtemp by 0.663in
    \rlap{\kern 0.100in\lower\graphtemp\hbox to 0pt{\hss $a$\hss}}%
    \graphtemp=.5ex\advance\graphtemp by 0.663in
    \rlap{\kern 0.350in\lower\graphtemp\hbox to 0pt{\hss $b$\hss}}%
    \graphtemp=.5ex\advance\graphtemp by 0.663in
    \rlap{\kern 0.600in\lower\graphtemp\hbox to 0pt{\hss $c$\hss}}%
    \graphtemp=.5ex\advance\graphtemp by 0.663in
    \rlap{\kern 0.850in\lower\graphtemp\hbox to 0pt{\hss $d$\hss}}%
    \graphtemp=.5ex\advance\graphtemp by 0.663in
    \rlap{\kern 1.100in\lower\graphtemp\hbox to 0pt{\hss $e$\hss}}%
    \graphtemp=.5ex\advance\graphtemp by 0.913in
    \rlap{\kern 0.600in\lower\graphtemp\hbox to 0pt{\hss $((abc)de) \to abc|de|$\hss}}%
    \graphtemp=.5ex\advance\graphtemp by 0.563in
    \rlap{\kern 1.600in\lower\graphtemp\hbox to 0pt{\hss $\bullet$\hss}}%
    \graphtemp=.5ex\advance\graphtemp by 0.563in
    \rlap{\kern 1.850in\lower\graphtemp\hbox to 0pt{\hss $\bullet$\hss}}%
    \graphtemp=.5ex\advance\graphtemp by 0.563in
    \rlap{\kern 2.100in\lower\graphtemp\hbox to 0pt{\hss $\bullet$\hss}}%
    \graphtemp=.5ex\advance\graphtemp by 0.563in
    \rlap{\kern 2.350in\lower\graphtemp\hbox to 0pt{\hss $\bullet$\hss}}%
    \graphtemp=.5ex\advance\graphtemp by 0.563in
    \rlap{\kern 2.600in\lower\graphtemp\hbox to 0pt{\hss $\bullet$\hss}}%
    \graphtemp=.5ex\advance\graphtemp by 0.063in
    \rlap{\kern 2.100in\lower\graphtemp\hbox to 0pt{\hss $\bullet$\hss}}%
    \graphtemp=.5ex\advance\graphtemp by 0.313in
    \rlap{\kern 2.100in\lower\graphtemp\hbox to 0pt{\hss $\bullet$\hss}}%
    \special{pa 1850 563}%
    \special{pa 2100 313}%
    \special{pa 2100 563}%
    \special{pa 2100 313}%
    \special{pa 2350 563}%
    \special{fp}%
    \special{pa 2100 313}%
    \special{pa 2100 63}%
    \special{pa 1600 563}%
    \special{pa 2100 63}%
    \special{pa 2600 563}%
    \special{fp}%
    \graphtemp=.5ex\advance\graphtemp by 0.663in
    \rlap{\kern 1.600in\lower\graphtemp\hbox to 0pt{\hss $a$\hss}}%
    \graphtemp=.5ex\advance\graphtemp by 0.663in
    \rlap{\kern 1.850in\lower\graphtemp\hbox to 0pt{\hss $b$\hss}}%
    \graphtemp=.5ex\advance\graphtemp by 0.663in
    \rlap{\kern 2.100in\lower\graphtemp\hbox to 0pt{\hss $c$\hss}}%
    \graphtemp=.5ex\advance\graphtemp by 0.663in
    \rlap{\kern 2.350in\lower\graphtemp\hbox to 0pt{\hss $d$\hss}}%
    \graphtemp=.5ex\advance\graphtemp by 0.663in
    \rlap{\kern 2.600in\lower\graphtemp\hbox to 0pt{\hss $e$\hss}}%
    \graphtemp=.5ex\advance\graphtemp by 0.913in
    \rlap{\kern 2.100in\lower\graphtemp\hbox to 0pt{\hss $(a(bcd)e) \to abcd|e|$\hss}}%
    \graphtemp=.5ex\advance\graphtemp by 0.563in
    \rlap{\kern 3.100in\lower\graphtemp\hbox to 0pt{\hss $\bullet$\hss}}%
    \graphtemp=.5ex\advance\graphtemp by 0.563in
    \rlap{\kern 3.350in\lower\graphtemp\hbox to 0pt{\hss $\bullet$\hss}}%
    \graphtemp=.5ex\advance\graphtemp by 0.563in
    \rlap{\kern 3.600in\lower\graphtemp\hbox to 0pt{\hss $\bullet$\hss}}%
    \graphtemp=.5ex\advance\graphtemp by 0.563in
    \rlap{\kern 3.850in\lower\graphtemp\hbox to 0pt{\hss $\bullet$\hss}}%
    \graphtemp=.5ex\advance\graphtemp by 0.563in
    \rlap{\kern 4.100in\lower\graphtemp\hbox to 0pt{\hss $\bullet$\hss}}%
    \graphtemp=.5ex\advance\graphtemp by 0.063in
    \rlap{\kern 3.350in\lower\graphtemp\hbox to 0pt{\hss $\bullet$\hss}}%
    \graphtemp=.5ex\advance\graphtemp by 0.313in
    \rlap{\kern 3.850in\lower\graphtemp\hbox to 0pt{\hss $\bullet$\hss}}%
    \special{pa 3600 563}%
    \special{pa 3850 313}%
    \special{pa 3850 563}%
    \special{pa 3850 313}%
    \special{pa 4100 563}%
    \special{fp}%
    \special{pa 3850 313}%
    \special{pa 3350 63}%
    \special{pa 3100 563}%
    \special{pa 3350 63}%
    \special{pa 3350 563}%
    \special{fp}%
    \graphtemp=.5ex\advance\graphtemp by 0.663in
    \rlap{\kern 3.100in\lower\graphtemp\hbox to 0pt{\hss $a$\hss}}%
    \graphtemp=.5ex\advance\graphtemp by 0.663in
    \rlap{\kern 3.350in\lower\graphtemp\hbox to 0pt{\hss $b$\hss}}%
    \graphtemp=.5ex\advance\graphtemp by 0.663in
    \rlap{\kern 3.600in\lower\graphtemp\hbox to 0pt{\hss $c$\hss}}%
    \graphtemp=.5ex\advance\graphtemp by 0.663in
    \rlap{\kern 3.850in\lower\graphtemp\hbox to 0pt{\hss $d$\hss}}%
    \graphtemp=.5ex\advance\graphtemp by 0.663in
    \rlap{\kern 4.100in\lower\graphtemp\hbox to 0pt{\hss $e$\hss}}%
    \graphtemp=.5ex\advance\graphtemp by 0.913in
    \rlap{\kern 3.600in\lower\graphtemp\hbox to 0pt{\hss $(ab(cde)) \to abcde||$\hss}}%
    \hbox{\vrule depth0.913in width0pt height 0pt}%
    \kern 4.200in
  }%
}%
}
\sp
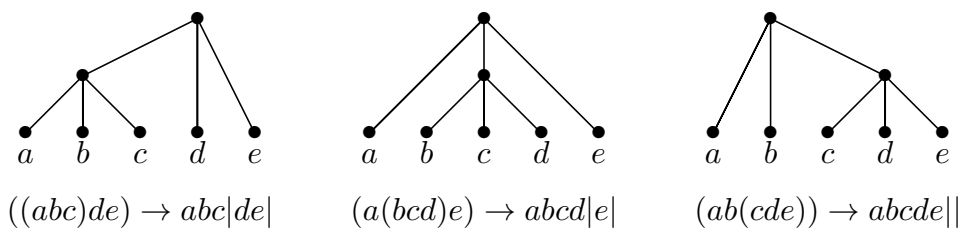
\ce{Figure 3.  The rooted plane ternary trees with five leaves.}

By working up from the leaves, combining $q+1$ subtrees at each non-leaf
vertex, we see that each tree corresponds to a product of $qn+1$ elements
in order using a non-associative $q+1$-ary operator.  Conversely, from such
a product we obtain the corresponding tree.  Sands [14] was interested in
counting the ways to form such a product.

In order to count the products, we convert the tree to a sequence of objects and
markers.  We do this by traversing the tree; beginning at the root, walk around
the tree, keeping our left hand on it, until we have traversed every edge twice.
We record an object for each leaf when it is visited and a marker for each
non-leaf vertex when the subtree rooted at that node is completed.  Figure 3
shows the resulting sequences when $q=n=2$.  The sequence corresponding to each
tree (that is, each bracketing) is ``post-fix'' notation for the formation of
the product; the marker specifies application of the operation to the $q$ most
recent objects on the stack, replacing them by a single object.

By converting markers to 1's and objects to 0's, we obtain an
$(n,qn+1)$-sequence.  For each such sequence produced in this way,
the number of objects that precede the $i$th marker must exceed $qi$,
since applying a $q+1$-ary operator $i$ times converts $qi+1$ objects into one
object.  Sands observed by an easy induction on $n$ that the $q$-dominating
condition is also sufficient.  Thus these bracketings (or these trees) are
equinumerous with the $q$-satisfying $(n,qn)$-sequences, and there are $C_n^q$
of them.

When $q=1$, we obtain the ballot sequences and the ordinary Catalan numbers.
Other counting problems solved by the Catalan numbers generalize in analogous
ways.  Hilton and Pedersen [8] observed that $C_n^q$ counts the subdivisions of
a convex polygon into $n$ disjoint $(q+1)$-gons by noncrossing diagonals.  This
generalizes a bijection between binary trees with $n+1$ leaves and dissections
of an $n+2$-gon into triangles, where the root becomes one edge of the polygon
and the leaves become the other edges.  In light of these generalizations,
we use the name $q$-{\it ballot sequences} for the $q$-satisfying sequences
with exactly $n$ 1's and $qn$ 0's.

\SH
{3. {\bigmit q}-SATISFYING SEQUENCES AND THE BRICKLAYER PROBLEM}
Before attacking the full generality of the bricklayer problem, we consider
how the bijection works when $q=1$.  Propp (see [12]) observed that the
ordinary Catalan numbers solve a simple coin-stacking problem.  We begin
with a base row of $n$ coins.  Each coin not in the base row rests on two coins
in the row immediately below it.  Figure 4 illustrates such a stack with a
base row of 6 coins.
\gpic{
\expandafter\ifx\csname graph\endcsname\relax \csname newbox\endcsname\graph\fi
\expandafter\ifx\csname graphtemp\endcsname\relax \csname newdimen\endcsname\graphtemp\fi
\setbox\graph=\vtop{\vskip 0pt\hbox{%
    \special{pn 8}%
    \special{pa 0 234}%
    \special{pa 281 234}%
    \special{fp}%
    \special{pa 312 234}%
    \special{pa 593 234}%
    \special{fp}%
    \special{pa 625 234}%
    \special{pa 906 234}%
    \special{fp}%
    \special{pa 937 234}%
    \special{pa 1218 234}%
    \special{fp}%
    \special{pa 1249 234}%
    \special{pa 1530 234}%
    \special{fp}%
    \special{pa 1561 234}%
    \special{pa 1842 234}%
    \special{fp}%
    \special{pa 312 78}%
    \special{pa 593 78}%
    \special{fp}%
    \special{pa 156 156}%
    \special{pa 437 156}%
    \special{fp}%
    \special{pa 468 156}%
    \special{pa 749 156}%
    \special{fp}%
    \special{pa 781 156}%
    \special{pa 1062 156}%
    \special{fp}%
    \special{pa 1405 156}%
    \special{pa 1686 156}%
    \special{fp}%
    \special{pa 2326 234}%
    \special{pa 2483 156}%
    \special{pa 2639 234}%
    \special{pa 2795 156}%
    \special{pa 2951 234}%
    \special{pa 3107 156}%
    \special{pa 3263 234}%
    \special{pa 3419 156}%
    \special{pa 3575 234}%
    \special{fp}%
    \special{pa 3575 234}%
    \special{pa 3732 156}%
    \special{pa 3888 234}%
    \special{pa 4044 156}%
    \special{pa 4200 234}%
    \special{fp}%
    \special{pa 2639 78}%
    \special{pa 2795 0}%
    \special{pa 2951 78}%
    \special{fp}%
    \special{pa 2483 156}%
    \special{pa 2639 78}%
    \special{pa 2795 156}%
    \special{pa 2951 78}%
    \special{pa 3107 156}%
    \special{pa 3263 78}%
    \special{pa 3419 156}%
    \special{fp}%
    \special{pa 3732 156}%
    \special{pa 3888 78}%
    \special{pa 4044 156}%
    \special{fp}%
    \hbox{\vrule depth0.234in width0pt height 0pt}%
    \kern 4.200in
  }%
}%
}
\sp
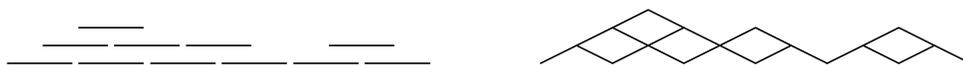
\ce{Figure 4.  A stack of coins and its conversion.}

Let $a_n$ be the number of distinct stacks that can be built on a
base of length $n$, with $a_0=1$.  If the number of contiguous coins at
the beginning of the first row above the base is $k-1$, then the stack is
completed by building one stack based on these $k-1$ coins and another based on
the last $n-k$ coins of the original base.  Summing over the possible values of
$k$ yields $a_n = \SE k1n a_{k-1} a_{n-k}$ for $n \ge 1$.  This is the
well-known recurrence satisfied by the Catalan numbers.  It arises for
the ballot sequence model by letting $2k$ be the minimum length of
a nonempty prefix with the same number of 0's and 1's.

The recurrence suggests a natural bijection.  Replace each coin by a wedge, as
shown on the right in Figure 4.  View each wedge as having width 2.  The wedge
for a coin above the base rests on the apexes of the wedges for its supporting
coins.  Now follow the outline of the mound of wedges.  Each step up is a 0,
each step down is a 1.  Each step moves one unit right, so there are $2n$ steps.
We end at the base, so there are $n$ steps of each type.  Since the mound never
dips below the base, the result is a ballot sequence.  That this is a bijection
is a special case of our main theorem.

We generalize to the bricklayer problem by viewing the coins as bricks of length
$2=q+1$.  We drop the requirement that the base be full to allow bases of all
lengths rather than merely multiples of $q+1$.  For this more general problem,
the proof is simpler.  In the special case where $m=(q+1)n$ and the base is
filled by $n$ bricks, there will be $C_n^q$ stacks, corresponding to the
$q$-ballot sequences.

Recall our conditions for $q$-stacks in the bricklayer problem:
\ti
1) The base has length $m$.
\ti
2) The bricks have length $q+1$.
\ti
3) Each brick not directly resting on the base rests on two contiguous bricks
\ti
immediately below it and covers a positive integer amount of the tops of each.
\br
Thus a brick has $q$ possible positions in relation to the two bricks below it.

\TH{3}
In the bricklayer problem with a base of length $m$, the number of $q$-stacks
with $n$ bricks resting directly on the base, and the total number of
$q$-stacks, respectively, are
$$\FR{m-(q+1)n+1}{m-n+1} \CH mn \qquad{\rm and}\qquad
\SE n0{\FL{m(q+1)}} \FR{m-(q+1)n+1}{m-n+1} \CH mn.$$
\PF
Fix $q$.  Let the {\it length} of a $q$-stack be the length of the base.
Let $S_m$ be the set of $q$-stacks of length $m$.
Let $Q_m$ be the set of $q$-satisfying sequences of length $m$.
We establish a bijection $\MAP f{S_m}{Q_m}$.
Furthermore, $f$ restricts to a bijection from $S_m^n$ to $Q_m^n$, where $S_m^n$
is the subset of $S_m$ consisting of stacks having $n$ bricks resting on the
base and $Q_m^n$ is the subset of $Q_m$ consisting of the sequences with $n$
1's.  By Corollary 2, establishing this bijection completes the proof.

To define $f$ on a stack $A\in S_m$, we begin by shaving the top corners of each
brick along a line from the bottom corner to a point at distance one from the
top corner along the top edge.  Each brick becomes a symmetric trapezoid with
upper edge of length $q-1$ and lower edge of length $q+1$ (a wedge when $q=1$).
Shaving the bricks is an invertible process (mathematically, not physically), so
we henceforth treat $S_m$ in this form.

For $A\in S_m$ (shaved), we define $f(A)$ by reading the top outline of the
stack, recording a 0 for each up-slant or horizontal step, and recording a 1 for
each down-slant.  Since we record one bit for each increase in the horizontal
coordinate, the length of $f(A)$ is $m$.  Figure 5 illustrates the
correspondence for one stack with $(q,m,n)=(2,12,3)$ and for all stacks with
$(q,m,n)=(2,9,3)$.
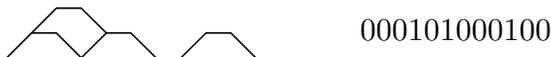
\gpic{
\expandafter\ifx\csname graph\endcsname\relax \csname newbox\endcsname\graph\fi
\expandafter\ifx\csname graphtemp\endcsname\relax \csname newdimen\endcsname\graphtemp\fi
\setbox\graph=\vtop{\vskip 0pt\hbox{%
    \special{pn 8}%
    \special{pa 0 217}%
    \special{pa 109 109}%
    \special{pa 217 109}%
    \special{pa 326 217}%
    \special{pa 435 109}%
    \special{pa 543 109}%
    \special{pa 652 217}%
    \special{pa 761 217}%
    \special{pa 870 109}%
    \special{fp}%
    \special{pa 870 109}%
    \special{pa 978 109}%
    \special{pa 1087 217}%
    \special{pa 1196 217}%
    \special{pa 1304 217}%
    \special{fp}%
    \special{pa 109 109}%
    \special{pa 217 0}%
    \special{pa 326 0}%
    \special{pa 435 109}%
    \special{fp}%
    \graphtemp=.5ex\advance\graphtemp by 0.109in
    \rlap{\kern 1.957in\lower\graphtemp\hbox to 0pt{\hss $000101000100$\hss}}%
    \hbox{\vrule depth0.217in width0pt height 0pt}%
    \kern 2.000in
  }%
}%
}

\gpic{
\expandafter\ifx\csname graph\endcsname\relax \csname newbox\endcsname\graph\fi
\expandafter\ifx\csname graphtemp\endcsname\relax \csname newdimen\endcsname\graphtemp\fi
\setbox\graph=\vtop{\vskip 0pt\hbox{%
    \special{pn 8}%
    \special{pa 0 255}%
    \special{pa 127 127}%
    \special{pa 255 127}%
    \special{pa 382 255}%
    \special{pa 510 127}%
    \special{pa 637 127}%
    \special{pa 764 255}%
    \special{pa 892 127}%
    \special{pa 1019 127}%
    \special{fp}%
    \special{pa 1019 127}%
    \special{pa 1146 255}%
    \special{fp}%
    \graphtemp=.5ex\advance\graphtemp by 0.127in
    \rlap{\kern 1.656in\lower\graphtemp\hbox to 0pt{\hss $001001001$\hss}}%
    \special{pa 0 764}%
    \special{pa 127 637}%
    \special{pa 255 637}%
    \special{pa 382 764}%
    \special{pa 510 637}%
    \special{pa 637 637}%
    \special{pa 764 764}%
    \special{pa 892 637}%
    \special{pa 1019 637}%
    \special{fp}%
    \special{pa 1019 637}%
    \special{pa 1146 764}%
    \special{fp}%
    \special{pa 127 637}%
    \special{pa 255 510}%
    \special{pa 382 510}%
    \special{pa 510 637}%
    \special{fp}%
    \graphtemp=.5ex\advance\graphtemp by 0.637in
    \rlap{\kern 1.656in\lower\graphtemp\hbox to 0pt{\hss $000101001$\hss}}%
    \special{pa 0 1274}%
    \special{pa 127 1146}%
    \special{pa 255 1146}%
    \special{pa 382 1274}%
    \special{pa 510 1146}%
    \special{pa 637 1146}%
    \special{pa 764 1274}%
    \special{pa 892 1146}%
    \special{pa 1019 1146}%
    \special{fp}%
    \special{pa 1019 1146}%
    \special{pa 1146 1274}%
    \special{fp}%
    \special{pa 255 1146}%
    \special{pa 382 1019}%
    \special{pa 510 1019}%
    \special{pa 637 1146}%
    \special{fp}%
    \graphtemp=.5ex\advance\graphtemp by 1.146in
    \rlap{\kern 1.656in\lower\graphtemp\hbox to 0pt{\hss $000011001$\hss}}%
    \special{pa 0 1783}%
    \special{pa 127 1656}%
    \special{pa 255 1656}%
    \special{pa 382 1783}%
    \special{pa 510 1656}%
    \special{pa 637 1656}%
    \special{pa 764 1783}%
    \special{pa 892 1656}%
    \special{pa 1019 1656}%
    \special{fp}%
    \special{pa 1019 1656}%
    \special{pa 1146 1783}%
    \special{fp}%
    \special{pa 510 1656}%
    \special{pa 637 1529}%
    \special{pa 764 1529}%
    \special{pa 892 1656}%
    \special{fp}%
    \graphtemp=.5ex\advance\graphtemp by 1.656in
    \rlap{\kern 1.656in\lower\graphtemp\hbox to 0pt{\hss $001000101$\hss}}%
    \special{pa 0 2293}%
    \special{pa 127 2166}%
    \special{pa 255 2166}%
    \special{pa 382 2293}%
    \special{pa 510 2166}%
    \special{pa 637 2166}%
    \special{pa 764 2293}%
    \special{pa 892 2166}%
    \special{pa 1019 2166}%
    \special{fp}%
    \special{pa 1019 2166}%
    \special{pa 1146 2293}%
    \special{fp}%
    \special{pa 637 2166}%
    \special{pa 764 2038}%
    \special{pa 892 2038}%
    \special{pa 1019 2166}%
    \special{fp}%
    \graphtemp=.5ex\advance\graphtemp by 2.166in
    \rlap{\kern 1.656in\lower\graphtemp\hbox to 0pt{\hss $001000011$\hss}}%
    \special{pa 0 2803}%
    \special{pa 127 2675}%
    \special{pa 255 2675}%
    \special{pa 382 2803}%
    \special{pa 510 2675}%
    \special{pa 637 2675}%
    \special{pa 764 2803}%
    \special{pa 892 2675}%
    \special{pa 1019 2675}%
    \special{fp}%
    \special{pa 1019 2675}%
    \special{pa 1146 2803}%
    \special{fp}%
    \special{pa 127 2675}%
    \special{pa 255 2548}%
    \special{pa 382 2548}%
    \special{pa 510 2675}%
    \special{pa 637 2548}%
    \special{pa 764 2548}%
    \special{pa 892 2675}%
    \special{fp}%
    \graphtemp=.5ex\advance\graphtemp by 2.675in
    \rlap{\kern 1.656in\lower\graphtemp\hbox to 0pt{\hss $000100101$\hss}}%
    \special{pa 2293 255}%
    \special{pa 2420 127}%
    \special{pa 2548 127}%
    \special{pa 2675 255}%
    \special{pa 2803 127}%
    \special{pa 2930 127}%
    \special{pa 3057 255}%
    \special{pa 3185 127}%
    \special{pa 3312 127}%
    \special{fp}%
    \special{pa 3312 127}%
    \special{pa 3439 255}%
    \special{fp}%
    \special{pa 2420 127}%
    \special{pa 2548 0}%
    \special{pa 2675 0}%
    \special{pa 2803 127}%
    \special{fp}%
    \special{pa 2930 127}%
    \special{pa 3057 0}%
    \special{pa 3185 0}%
    \special{pa 3312 127}%
    \special{fp}%
    \graphtemp=.5ex\advance\graphtemp by 0.127in
    \rlap{\kern 3.949in\lower\graphtemp\hbox to 0pt{\hss $000100011$\hss}}%
    \special{pa 2293 764}%
    \special{pa 2420 637}%
    \special{pa 2548 637}%
    \special{pa 2675 764}%
    \special{pa 2803 637}%
    \special{pa 2930 637}%
    \special{pa 3057 764}%
    \special{pa 3185 637}%
    \special{pa 3312 637}%
    \special{fp}%
    \special{pa 3312 637}%
    \special{pa 3439 764}%
    \special{fp}%
    \special{pa 2548 637}%
    \special{pa 2675 510}%
    \special{pa 2803 510}%
    \special{pa 2930 637}%
    \special{pa 3057 510}%
    \special{pa 3185 510}%
    \special{pa 3312 637}%
    \special{fp}%
    \graphtemp=.5ex\advance\graphtemp by 0.637in
    \rlap{\kern 3.949in\lower\graphtemp\hbox to 0pt{\hss $001001001$\hss}}%
    \special{pa 2293 1274}%
    \special{pa 2420 1146}%
    \special{pa 2548 1146}%
    \special{pa 2675 1274}%
    \special{pa 2803 1146}%
    \special{pa 2930 1146}%
    \special{pa 3057 1274}%
    \special{pa 3185 1146}%
    \special{pa 3312 1146}%
    \special{fp}%
    \special{pa 3312 1146}%
    \special{pa 3439 1274}%
    \special{fp}%
    \special{pa 2420 1146}%
    \special{pa 2548 1019}%
    \special{pa 2675 1019}%
    \special{pa 2803 1146}%
    \special{pa 2930 1019}%
    \special{pa 3057 1019}%
    \special{pa 3185 1146}%
    \special{fp}%
    \special{pa 2548 1019}%
    \special{pa 2675 892}%
    \special{pa 2803 892}%
    \special{pa 2930 1019}%
    \special{fp}%
    \graphtemp=.5ex\advance\graphtemp by 1.146in
    \rlap{\kern 3.949in\lower\graphtemp\hbox to 0pt{\hss $000010011$\hss}}%
    \special{pa 2293 1783}%
    \special{pa 2420 1656}%
    \special{pa 2548 1656}%
    \special{pa 2675 1783}%
    \special{pa 2803 1656}%
    \special{pa 2930 1656}%
    \special{pa 3057 1783}%
    \special{pa 3185 1656}%
    \special{pa 3312 1656}%
    \special{fp}%
    \special{pa 3312 1656}%
    \special{pa 3439 1783}%
    \special{fp}%
    \special{pa 2420 1656}%
    \special{pa 2548 1529}%
    \special{pa 2675 1529}%
    \special{pa 2803 1656}%
    \special{pa 2930 1529}%
    \special{pa 3057 1529}%
    \special{pa 3185 1656}%
    \special{fp}%
    \special{pa 2675 1529}%
    \special{pa 2803 1401}%
    \special{pa 2930 1401}%
    \special{pa 3057 1529}%
    \special{fp}%
    \graphtemp=.5ex\advance\graphtemp by 1.656in
    \rlap{\kern 3.949in\lower\graphtemp\hbox to 0pt{\hss $000010101$\hss}}%
    \special{pa 2293 2293}%
    \special{pa 2420 2166}%
    \special{pa 2548 2166}%
    \special{pa 2675 2293}%
    \special{pa 2803 2166}%
    \special{pa 2930 2166}%
    \special{pa 3057 2293}%
    \special{pa 3185 2166}%
    \special{pa 3312 2166}%
    \special{fp}%
    \special{pa 3312 2166}%
    \special{pa 3439 2293}%
    \special{fp}%
    \special{pa 2548 2166}%
    \special{pa 2675 2038}%
    \special{pa 2803 2038}%
    \special{pa 2930 2166}%
    \special{pa 3057 2038}%
    \special{pa 3185 2038}%
    \special{pa 3312 2166}%
    \special{fp}%
    \special{pa 2675 2038}%
    \special{pa 2803 1911}%
    \special{pa 2930 1911}%
    \special{pa 3057 2038}%
    \special{fp}%
    \graphtemp=.5ex\advance\graphtemp by 2.166in
    \rlap{\kern 3.949in\lower\graphtemp\hbox to 0pt{\hss $000001011$\hss}}%
    \special{pa 2293 2803}%
    \special{pa 2420 2675}%
    \special{pa 2548 2675}%
    \special{pa 2675 2803}%
    \special{pa 2803 2675}%
    \special{pa 2930 2675}%
    \special{pa 3057 2803}%
    \special{pa 3185 2675}%
    \special{pa 3312 2675}%
    \special{fp}%
    \special{pa 3312 2675}%
    \special{pa 3439 2803}%
    \special{fp}%
    \special{pa 2548 2675}%
    \special{pa 2675 2548}%
    \special{pa 2803 2548}%
    \special{pa 2930 2675}%
    \special{pa 3057 2548}%
    \special{pa 3185 2548}%
    \special{pa 3312 2675}%
    \special{fp}%
    \special{pa 2803 2548}%
    \special{pa 2930 2420}%
    \special{pa 3057 2420}%
    \special{pa 3185 2548}%
    \special{fp}%
    \graphtemp=.5ex\advance\graphtemp by 2.675in
    \rlap{\kern 3.949in\lower\graphtemp\hbox to 0pt{\hss $000000111$\hss}}%
    \hbox{\vrule depth2.803in width0pt height 0pt}%
    \kern 4.000in
  }%
}%
}
\sp
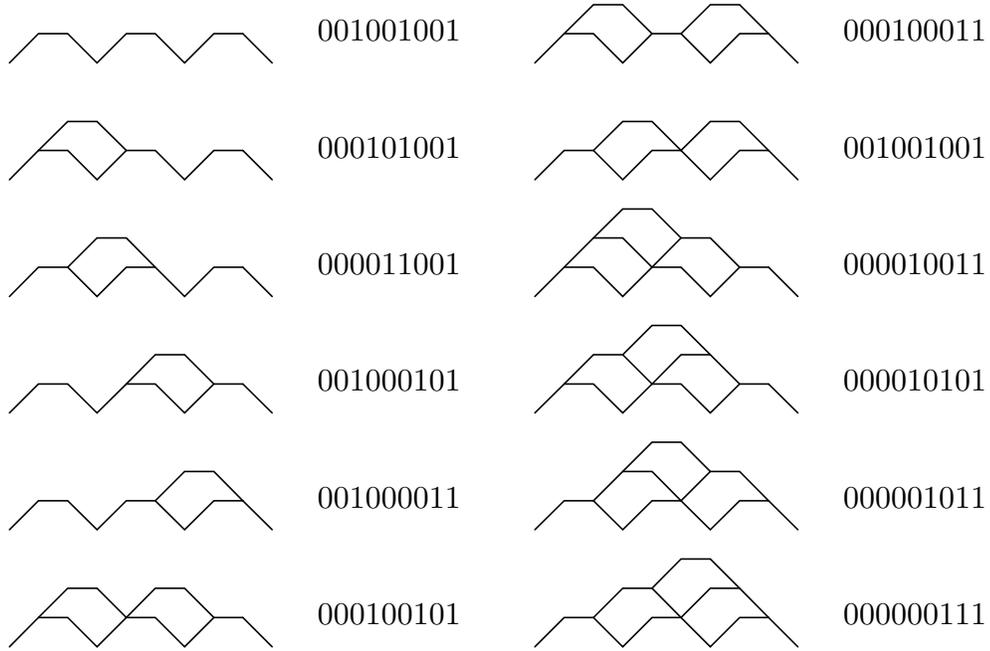
\ce{Figure 5.  Shaved $q$-stacks and the corresponding $q$-satisfying sequences.}

Among the $q$-satisfying sequences, the $q$-dominating sequences are those that
have no $q$-ballot prefix.  We prove the further property that, within each pair
$(S_m^n, Q_m^n)$, $f$ matches the $q$-dominating sequences with the stacks not
covering the first space on the base.  We use induction on $m$.  When $m=1$,
there are no bricks, and the empty stack in $S_1^0$ maps to the sequence 0 in
$Q_1^0$.  For $m>1$, we consider three types of stacks and the corresponding
sequences, showing first that $f$ restricts as desired.

If $A\in S_m^n$ does not cover the first space of the base, then $f(A)$ consists
of 0 followed by $f(A')$, where $A'$ is obtained from $A$ by deleting the first
space of the base.  By the induction hypothesis, $f(A')\in Q_{m-1}^n$, so
$f(A)\in Q_m^n$ and $f(A)$ is $q$-dominating.

If $A\in S_m^n$ covers the first space, let $m'=k(q+1)$ be the step on
which the outline of $A$ first returns to the base (in the top example of
Figure 5, $k=2$ and $m'=6$).  If $m'<m$, then $A$ is the concatenation of a
stack $A'\in S_{m'}^k$ and a stack $A''\in S_{m-m'}^{n-k}$, and $f(A)$ is the
concatenation of $f(A')$ and $f(A'')$.  By the induction hypothesis, $f(A')$ is
a $q$-ballot sequence and $f(A'')$ is $q$-satisfying, so $f(A)\in Q_m^n$ and
$f(A)$ is not $q$-dominating.

If $m'=m$, then $m=(q+1)n$ and every notch in the lowest row of bricks is
covered by a higher brick.  Let $A'$ be the stack obtained from $A$ by deleting
the bottom row of bricks and the first and last space of the base.  Since every
brick above the first row covers one notch, $A'\in S_{m-2}^{n-1}$.  Also, $f(A)$
is obtained from $f(A')$ by adding 0 at the beginning and 1 at the end.  By the
induction hypothesis, $f(A')\in Q_{m-2}^{n-1}$.  Adding 0 at the beginning and 1
at the end of a $q$-satisfying sequence with these parameters yields a
$q$-satisfying sequence, so $f(A)\in Q_{(q+1)n}^n$.  The prefix without the
final 1 is $q$-dominating and has no proper $q$-ballot prefix, but the full
sequence $f(A)$ is a $q$-ballot sequence and hence is not $q$-dominating.

To show that $f$ is invertible, consider $a\in Q_m^n$.  If $a$ is
$q$-dominating, then no stack covering the first space has image $a$.  Let $a'$
be the $q$-statisfying sequence obtained by deleting the initial 0 of $a$.  By
the induction hypothesis, there exists one stack $A'\in S_{m-1}^n$ such that
$f(A')=a'$.  Among the stacks in $S_m^n$ not covering the first space, there is
thus one stack $A$ such that $f(A)=a$.

If $a$ is not $q$-dominating, then $a$ has a $q$-ballot prefix.  Let $m'$ be the
length of the shortest $q$-ballot prefix $a'$ of $a$, having $k$ 1's.  Since
$a'$ has $qk$ 0's, the remainder $a''$ of $a$ is in $Q_{m-m'}^{n-k}$.
We have observed that a proper $q$-ballot prefix of length $m'<m$ arises in
$f(A)$ if and only the outline of $A$ first returns to the base after $m'$
steps.  By the induction hypothesis, there is one $A\in S_{m'}^k$ such that
$f(A)=a'$, and there is one $A''\in A_{m-m'}^{n-k}$ such that $f(A'')=a''$,
and the concatenation yields the unique $A$ such that $f(A)=a$.

If $a$ is not $q$-dominating and $m'=m$, then $m=(q+1)n$ and $a$ is a $q$-ballot
sequence with no proper $q$-ballot prefix.  Before the final 1, $a$ is
$q$-dominating, and deleting the initial 0 and final 1 yields a $q$-satisfying
sequence $a'$ of length $m-2$.  We have observed that if $f(A)$ is a $q$-ballot
sequence with no proper $q$-ballot prefix, then $A$ has $n-1$ bricks in the
second row, and $f(A)$ is obtained by adding 0 at the beginning and 1 at the end
of $f(A')$ for the stack $A'\in S_{m-2}^{n-1}$ obtained by deleting the bottom
row of $A$ and shortening the base at both ends.  By the induction hypothesis,
there is one $A'\in S_{m-2}^{n-1}$ such that $f(A')=a'$, and thus there is one
$A$ such that $f(A)=a$.  \qed

The proof yields a recurrence for the number $c_{m,n}$ of $q$-satisfying
sequences of length $m$ with $n$ 1's.  We have $c_{1,0}=1$, and $c_{1,n}=0$ for
$n\ne1$.  For $m>1$,
$$c_{m,n} = c_{m-1,n} + \eps_{m,n}c_{m-2,n-1}
+\sum_{0\le k<m/(q+1)} c_{(q+1)k,k} c_{m-(q+1)k,n-k},\eqno{(*)}$$
where $\eps_{m,n}$ is 1 if $m=(q+1)n$ and is 0 otherwise.
Bailey [1] obtained another recurrence for the case $q=1$ and used it to
obtain the first statement of Corollary 2 in that case.

The Catalan recurrence $a_n = \SE k1n a_{k-1}a_{n-k}$ is simpler than $(*)$
because removing the initial 0 and trailing 1 from a ballot sequence that
has no balanced prefix yields a shorter ballot sequence.  This statement does
not generalize; when a $q$-ballot sequence of length $(q+1)n$ has no proper
$q$-ballot prefix (such as 000100011 for $q=2$ and $n=3$), removing the bits
after the penultimate 1 and removing enough leading 0's to reduce to length
$(q+1)(n-1)$ need not produce a $q$-ballot sequence.

We can obtain a natural recurrence for generalized Catalan numbers by modeling
the formation of bracketings.  We know that $C_n^q$ counts both the $q$-ballot
sequences with $n$ 1's and the bracketings of a product involving $n$
applications of a $q+1$-ary operator.  Hilton and Pedersen [8] observed that the
last application of the operator combines $q+1$ segments, each of which is a
shorter $q$-dominating $(n_i,qn_i+1)$-sequence.  Thus 
$$C_n^q = \sum \PE i1{q+1} C_{n_i}^q,$$
where the summation runs over all choices of $q+1$ positive integers
$\VEC n1{q+1}$ that sum to $n-1$.  The initial condition is $C_0^q=1$.

\SH
{4. ({\bigmit k},{\bigmit qk}+1)-ARRANGEMENTS}
For $p=1$, Kierstead and Trotter strengthened the Cycle Lemma on
$(k,qk+p)$-arrangements, showing that each 0 plays a special role.
Within a cyclic arrangement $a=\VEC a0{n-1}$ of 0's and 1's, we denote the list
$\VEC a{i+1}j$ (indices modulo $n$) by $(i,j]$.  The {\it linearization} of $a$
ending at position $i$ is $(i,i]$.   For a fixed linearization, we use {\it
0-interval} to mean a prefix ending at a 0.  If $a_i = 0$, then $(i,i]$ is a
{\it 0-linearization}, and the full list $(i,i]$ is the {\it trivial}
0-interval.  We use $w_0(I)$ and $w_1(I)$ to denote the number of 0's and 1's
in $I$, respectively.  A 0-interval $I$ is {\it q-good} if $w_0(I)>qw_1(I)$.
In every 0-linearization of a $(k,qk+p)$-arrangement with $p>0$, the trivial
0-interval is $q$-good.

For $q=p=1$ and $1\le i\le k+1$, Kierstead and Trotter [9] proved that every
$(k,k+1)$-arrangement $a$ has a unique 0-linearization such that exactly $i$ of
the 0-intervals are 1-good.  They noted that this result is implicit in the work
of Feller [6] and Narayana [11], and they used it to construct new explicit
perfect matchings in the bipartite graph of the inclusion relation on the
$k$-sets and $k+1$-sets of a $2k+1$-element set.  Their elegant proof of a
technically stronger statement extends
directly to $(k,qk+1)$-arrangements.  Figure 6 shows a (clockwise)
$(3,7)$-arrangement and its 0-linearizations, indicating the number of $2$-good
0-intervals in each and underscoring the positions that end $2$-good
0-intervals.

\def\mmat{\matrix{
1&\ 1\ 0\ 1\ 0\ 0\ 1\ 0\ 0\ 0\ \un0\cr
 &\cr
3&\ 1\ 0\ 0\ 1\ 0\ 0\ \un0\ \un0\ 1\ \un0\cr
 &\cr
6&\ \un0\ 1\ 0\ \un0\ \un0\ \un0\ 1\ \un0\ 1\ \un0\cr
4&\ 1\ 0\ 0\ \un0\ \un0\ 1\ \un0\ 1\ 0\ \un0\cr
 &\cr
7&\ \un0\ \un0\ \un0\ 1\ \un0\ 1\ \un0\ \un0\ 1\ \un0\cr
5&\ \un0\ \un0\ 1\ \un0\ 1\ 0\ \un0\ 1\ 0\ \un0\cr
2&\ \un0\ 1\ 0\ 1\ 0\ 0\ 1\ 0\ 0\ \un0\cr } }
\gpic{
\expandafter\ifx\csname graph\endcsname\relax \csname newbox\endcsname\graph\fi
\expandafter\ifx\csname graphtemp\endcsname\relax \csname newdimen\endcsname\graphtemp\fi
\setbox\graph=\vtop{\vskip 0pt\hbox{%
    \graphtemp=.5ex\advance\graphtemp by 0.206in
    \rlap{\kern 0.697in\lower\graphtemp\hbox to 0pt{\hss 1\hss}}%
    \graphtemp=.5ex\advance\graphtemp by 0.305in
    \rlap{\kern 1.001in\lower\graphtemp\hbox to 0pt{\hss 0\hss}}%
    \graphtemp=.5ex\advance\graphtemp by 0.563in
    \rlap{\kern 1.188in\lower\graphtemp\hbox to 0pt{\hss 1\hss}}%
    \graphtemp=.5ex\advance\graphtemp by 0.882in
    \rlap{\kern 1.188in\lower\graphtemp\hbox to 0pt{\hss 0\hss}}%
    \graphtemp=.5ex\advance\graphtemp by 1.140in
    \rlap{\kern 1.001in\lower\graphtemp\hbox to 0pt{\hss 0\hss}}%
    \graphtemp=.5ex\advance\graphtemp by 1.239in
    \rlap{\kern 0.697in\lower\graphtemp\hbox to 0pt{\hss 1\hss}}%
    \graphtemp=.5ex\advance\graphtemp by 1.140in
    \rlap{\kern 0.394in\lower\graphtemp\hbox to 0pt{\hss 0\hss}}%
    \graphtemp=.5ex\advance\graphtemp by 0.882in
    \rlap{\kern 0.206in\lower\graphtemp\hbox to 0pt{\hss 0\hss}}%
    \graphtemp=.5ex\advance\graphtemp by 0.563in
    \rlap{\kern 0.206in\lower\graphtemp\hbox to 0pt{\hss 0\hss}}%
    \graphtemp=.5ex\advance\graphtemp by 0.305in
    \rlap{\kern 0.394in\lower\graphtemp\hbox to 0pt{\hss 0\hss}}%
    \graphtemp=.5ex\advance\graphtemp by 0.722in
    \rlap{\kern 2.245in\lower\graphtemp\hbox to 0pt{\hss $\mmat$\hss}}%
    \graphtemp=.5ex\advance\graphtemp by 0.722in
    \rlap{\kern 3.794in\lower\graphtemp\hbox to 0pt{\hss  \hss}}%
    \hbox{\vrule depth1.445in width0pt height 0pt}%
    \kern 4.000in
  }%
}%
}
\sp
\ce{Figure 6.  2-good 0-intervals in 0-linearizations of a (3,7)-arrangement.}

\LM{4}
(Strong Cycle Lemma)  If $a$ is a $(k,qk+1)$-arrangement and $1\le i\le qk+1$,
then there is a unique 0-linearization $(j_i,j_i]$ of
$a$ in which exactly $i$ 0-intervals are $q$-good.  Furthermore, for $i\le qk$,
the 0's that end $q$-good 0-intervals in $(j_i,j_i]$ also end $q$-good
0-intervals in $(j_{i+1},j_{i+1}]$ (we call this the {\it nesting property}).
\PF
Given a $(k,qk+1)$-arrangement $a$, let the {\it deficiency} of an interval
$I=(r,s]$ be $\dlt(r,s] = qw_1(I)-w_0(I)$.  Given a 0-linearization $(r,r]$, let
$D(r) = \SET j:{\dlt(r,j]<0 \AND a_j=0}$; this is the set of indices ending
$q$-good 0-intervals for $(r,r]$.  Note that $r\in D(r)$.  Since
$1\le \C{D(r)}\le qk+1$ for all $j$, it suffices to prove that the sets
$\{D(j)\}$ are distinct and are linearly ordered by inclusion.

Let $r,s$ be the positions of two 0's.  Since $\dlt (r,s] + \dlt (s,r] = -1$,
exactly one of $(r,s]$ and $(s,r]$ has negative deficiency and is $q$-good.
When $(r,s]$ is $q$-good, we claim that $D(s)\subset D(r)$.  Note that
$s\in D(r)$, but $r\notin D(s)$.  Now consider $j\in D(s)$; we have two cases.
If $j\in(r,s]$, then $\dlt(r,j] = \dlt(s,j]-\dlt(s,r] < 0$.  If $j\in(s,r]$,
then $\dlt(r,j] = \dlt(s,j]+\dlt(r,s] < 0$.  In each case, we obtain
$j\in D(r)$.  \qed 

The Strong Cycle Lemma can also be proved by constructing $j_{i+1}$ explicitly
from $j_i$, but that takes longer.
\skipit{We omit this, because our objective is to show the power of
the Strong Cycle Lemma and the nesting property in applications related to
generalizations of the Catalan numbers.}
As an application, we obtain a result of Chung and Feller that generalizes the
Ballot Problem discussed earlier.  We obtain the simple Ballot Problem by
setting $l=0$ and interchanging A and B.  The Chung-Feller proof used analytic
methods.

\CO{5}
(Chung-Feller [2])  Let $l$ be an integer in $\{0,\ldots,n\}$.  In a random
sequence of $n$ A's and $n$ B's, the probability is $1/(n+1)$ that there are
exactly $l$ choices of $i$ such that the $i$th A precedes the $i$th B.
\PF
Given a sequence $b$ of $n$ A's and $n$ B's, convert A's to 0's, B's to 1's,
and append a 0 at the end; call this sequence $b'$.  The sequence $b'$ is a
0-linearization of an $(n,n+1)$-arrangement $a$.  Because $n+1$ and $2n+1$
are relatively prime, exactly $n+1$ sequences of $n$ A's and $n$ B's yield
the same $(n,n+1)$-arrangement.  By the Strong Cycle Lemma, the $n+1$
0-linearizations of $a$ have different numbers of 1-good 0-intervals.
The $i$th B in $b$ precedes the $i$th A in $b$ if and only if the $i$th 0 in
$b'$ is a 1-good 0-interval.  By grouping the sequences into sets of size
$n+1$ yielding the same cyclic arrangement, we see that the number of 
sequences with exactly $l$ values where the $i$th A precedes the $i$th B
is independent of $l$.  \qed

We next extend the Strong Cycle Lemma by specifying an arbitrary set of 0's
in a $(k,qk+1)$-arrangement.

\LM{6}
(Stronger Cycle Lemma)  If $a$ is a $(k,qk+1)$-arrangement, $S$ is a set of
$t$ positions containing 0's, and $1\le i\le t$, then there is a unique
0-linearization $(j_i,j_i]$ of $a$ such that $j_i\in S$ and exactly $i$ of the
0-intervals ending at elements of $S$ are $q$-good.  Furthermore, for $i < t$,
the 0's that end $q$-good 0-intervals in $(j_i,j_i]$ also end $q$-good
0-intervals in $(j_{i+1},j_{i+1}]$.
\PF
Define deficiency as in Lemma 4, but let $D(r)=\SET j:{\dlt(r,j]<0\AND j\in S}$.
Since $1\le  D(r)\le t$, it suffices to show that the sets $D(r)$ for $r\in S$
are distinct and ordered by inclusion.  Given $r,s \in S$, the proof of this
is as in Lemma 4.  \qed 

Viewing the elements of a $(k,k+1)$-arrangement as exponents on $-1$ yields a
cyclic arrangement of $k+1$ positive $1$'s and $k$ negative $1$'s.  The ordinary
Cycle Lemma provides a unique starting position such that all the partial sums
are positive.  Raney [13] proved more generally that every cyclic arrangement of
integers summing to $+1$ has a unique starting position such that all partial
sums are positive.  The Stronger Cycle Lemma for $q = 1$ provides a short proof
of a further generalization.  In Figure 7, we indicate the number of
positive partial sums in each successive linearization of a clockwise
arrangement and underscore the positions that end positive partial sums.

\def\mmat{\matrix{
5&&\ \un2&\un{-1}&\ \un2&-5&\ \un3&-2&\ 1&-2&\ \un3\cr
2&& -1 &\ \un2   &-5&\ 3&-2&\ 1&-2&\ 3&\ \un2\cr
3&&\ \un2&  -5   &\ 3&-2&\ 1&-2&\ 3&\ \un2&\un{-1}\cr
1&& -5 & \   3   &-2&\ 1&-2&\ 3&\ 2&{-1}&\ \un2\cr
8&&\ \un3&\un{-2}&\ \un1&-2&\ \un3&\ \un2&\un{-1}&\ \un2&\un{-5}\cr
4&& -2 &  \  1   &-2&\ 3&\ \un2&\un{-1}&\ \un2&-5&\ \un3\cr
7&&\ \un1&  -2   &\ \un3&\ \un2&\un{-1}&\ \un2&{-5}&\ \un3&\un{-2}\cr
6&& -2 &\ \un3   &\ \un2&\un{-1}&\ \un2&{-5}&\ \un3&{-2}&\ \un1\cr
9&&\ \un3&\ \un2 &\un{-1}&\ \un2&\un{-5}&\ \un3&\un{-2}&\ \un1&\un{-2}\cr
}  }
\gpic{
\expandafter\ifx\csname graph\endcsname\relax \csname newbox\endcsname\graph\fi
\expandafter\ifx\csname graphtemp\endcsname\relax \csname newdimen\endcsname\graphtemp\fi
\setbox\graph=\vtop{\vskip 0pt\hbox{%
    \graphtemp=.5ex\advance\graphtemp by 0.164in
    \rlap{\kern 0.566in\lower\graphtemp\hbox to 0pt{\hss 2\hss}}%
    \graphtemp=.5ex\advance\graphtemp by 0.259in
    \rlap{\kern 0.829in\lower\graphtemp\hbox to 0pt{\hss $-1$\hss}}%
    \graphtemp=.5ex\advance\graphtemp by 0.501in
    \rlap{\kern 0.969in\lower\graphtemp\hbox to 0pt{\hss 2\hss}}%
    \graphtemp=.5ex\advance\graphtemp by 0.777in
    \rlap{\kern 0.920in\lower\graphtemp\hbox to 0pt{\hss $-5$\hss}}%
    \graphtemp=.5ex\advance\graphtemp by 0.957in
    \rlap{\kern 0.706in\lower\graphtemp\hbox to 0pt{\hss 3\hss}}%
    \graphtemp=.5ex\advance\graphtemp by 0.957in
    \rlap{\kern 0.426in\lower\graphtemp\hbox to 0pt{\hss $-2$\hss}}%
    \graphtemp=.5ex\advance\graphtemp by 0.777in
    \rlap{\kern 0.212in\lower\graphtemp\hbox to 0pt{\hss 1\hss}}%
    \graphtemp=.5ex\advance\graphtemp by 0.501in
    \rlap{\kern 0.164in\lower\graphtemp\hbox to 0pt{\hss $-2$\hss}}%
    \graphtemp=.5ex\advance\graphtemp by 0.259in
    \rlap{\kern 0.303in\lower\graphtemp\hbox to 0pt{\hss 3\hss}}%
    \graphtemp=.5ex\advance\graphtemp by 0.572in
    \rlap{\kern 2.815in\lower\graphtemp\hbox to 0pt{\hss $\mmat$\hss}}%
    \graphtemp=.5ex\advance\graphtemp by 0.572in
    \rlap{\kern 3.836in\lower\graphtemp\hbox to 0pt{\hss  \hss}}%
    \hbox{\vrule depth1.120in width0pt height 0pt}%
    \kern 4.000in
  }%
}%
}
\sp
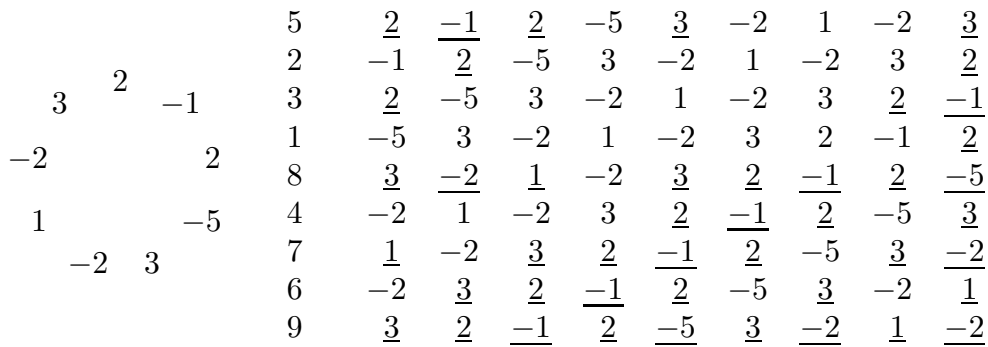
\ce{Figure 7.  Positive partial sums in an arrangement summing to $+1$.}

\CO{7}
(Mont\'agh [10])
Given a cyclic arrangement of $n$ integers summing to $+1$ and an integer
$l \in \{1,\ldots,n\}$, there is a unique linearization of the arrangement
such that exactly $l$ of the partial sums are positive.
\PF
Let $b$ denote the arrangement of integers.  Form a cyclic arrangement $a$ of
1's and 0's by replacing each nonnegative integer $b_i$ by a single 1 followed
by $1+b_i$ consecutive 0's, and replacing each negative $b_i$ by $1-b_i$ 
consecutive 1's followed by one 0.  The resulting $a$ is a
$(k,k+1)$-arrangement, where $k=n+\FL{\SG\C{b_i}/2}$.  Let $S$ be the $n$-set of
0's that end maximal consecutive segments of 0's.  The 0-linearizations ending
in $S$ correspond naturally to linearizations of $b$; the number of positive
partial sums in a linearization of $b$ equals the number of $1$-good 0-intervals
ending in $S$ in the corresponding 0-linearization of $a$.  By Lemma 6, these
numbers are distinct for the $n$ linearizations ending in $S$.  \qed 

In this application, the ``nesting property'' says that the numbers ending
positive partial sums in the $l-1$th arrangement also end positive partial sums
in the $l$th arrangement.

Graham, Knuth, and Patashnik [7, p.~346] presented a geometric proof of Raney's
original result, which upon closer examination also yields Mont\'agh's
generalization.  (Dershowitz and Zaks [3] observed that the geometric approach
can also be used to prove the Cycle Lemma itself.)  Encode the integer
arrangement $\VEC a1n$ as a walk in the plane, starting from the origin and
moving $(+1,+a_i)$ from the current position when the $i$th number is
encountered.  The ending position is $(n,1)$.  Figure 8 shows two periods of the
walk for the sequence $2,-1,2,-5,3,-2,1,-2,3$.

\gpic{
\expandafter\ifx\csname graph\endcsname\relax \csname newbox\endcsname\graph\fi
\expandafter\ifx\csname graphtemp\endcsname\relax \csname newdimen\endcsname\graphtemp\fi
\setbox\graph=\vtop{\vskip 0pt\hbox{%
    \graphtemp=.5ex\advance\graphtemp by 0.649in
    \rlap{\kern 0.032in\lower\graphtemp\hbox to 0pt{\hss $\bullet$\hss}}%
    \graphtemp=.5ex\advance\graphtemp by 0.519in
    \rlap{\kern 1.201in\lower\graphtemp\hbox to 0pt{\hss $\bullet$\hss}}%
    \graphtemp=.5ex\advance\graphtemp by 0.390in
    \rlap{\kern 2.370in\lower\graphtemp\hbox to 0pt{\hss $\bullet$\hss}}%
    \special{pn 8}%
    \special{pa 32 0}%
    \special{pa 32 649}%
    \special{pa 2500 649}%
    \special{fp}%
    \special{pa 1201 0}%
    \special{pa 1201 649}%
    \special{fp}%
    \special{pa 2370 0}%
    \special{pa 2370 649}%
    \special{fp}%
    \special{pa 32 649}%
    \special{pa 162 390}%
    \special{pa 292 519}%
    \special{pa 422 260}%
    \special{pa 552 909}%
    \special{pa 682 519}%
    \special{pa 812 779}%
    \special{pa 942 649}%
    \special{pa 1071 909}%
    \special{fp}%
    \special{pa 1071 909}%
    \special{pa 1201 519}%
    \special{pa 1331 260}%
    \special{pa 1461 390}%
    \special{pa 1591 130}%
    \special{pa 1721 779}%
    \special{pa 1851 390}%
    \special{pa 1981 649}%
    \special{pa 2110 519}%
    \special{fp}%
    \special{pa 2110 519}%
    \special{pa 2240 779}%
    \special{pa 2370 390}%
    \special{fp}%
    \hbox{\vrule depth0.909in width0pt height 0pt}%
    \kern 2.500in
  }%
}%
}
\sp
\ce{Figure 8.  A geometric argument.}

The unique starting position from which all the partial sums are positive
follows the last occurrence of the minimum in the first period.
All other positions have a non-positive partial sum ending at that position,
but partial sums starting after it are positive.  All but one partial sum is
positive when we start after the previous occurrence of the minimum or, when the
minimum is unique, after the last occurrence of the next smallest value.
For each $l<n$, let $b_l$ be the position in the first period from which $l+1$
partial sums are positive.  Then $b_{l-1}$ is obtained from $b_l$ by moving to
the previous occurrence of the same height as $b_l$ or, if $b_l$ is the first
occurrence of that height, the last occurrence of the next larger height.

When $n=(q+1)k+1$ and the sequence consists only of 1's and $-q$'s, summing to
$+1$ requires that exactly $k$ terms equal $-q$.  Following Graham, Knuth, and
Patashnik (with a shift of index), we call such a sequence a $q$-{\it Raney
sequence} if all the partial sums are positive (this requires the first term to
be a 1).  They prove that there are $C_k^q$ such sequences.  This follows
immediately from Corollary 7 when $l=n$, since $C_k^q$ equals the number of
$(k,qk+1)$-arrangements.

The Strong Cycle Lemma also yields a short direct proof that the number of
$q$-ballot sequences with $k$ ones is $C_k^q$.  As before, prepending a 0 shows
that these are equinumerous with the $q$-dominating $(k,qk+1)$-sequences.  By
the Strong Cycle Lemma, the reverse $a'$ of such a sequence $a$ has a unique
0-linearization such that {\it no} 0-interval is $q$-good.  The reverse of this
0-linearization is the unique cyclic permutation of $a$ such that
$w_0(I)>qw_1(I)$ for every prefix $I$ (whether ending at a 0 or a 1).
Hence we conclude again that the $q$-ballot sequences of length $(q+1)k$
are equinumerous with the $(k,qk+1)$-arrangements.

The result of Kierstead and Trotter, extended to the Strong Cycle Lemma, 
distinguishes the 0's of a $(k,qk+1)$-arrangement in a combinatorial fashion.
We close this section by presenting another combinatorial distinguishing of
these 0's that extends to $(k,l)$-arrangements whenever $k$ and $l$ are
relatively prime.
In the case $l=k+1$, it yields matchings different from the matchings of
Kierstead and Trotter [9] between the middle levels of the lattice of subsets of
a $2k+1$-element set; further discussion appears in [4].

\TH{8}
(Snevily [15])
If $k$ and $l$ are relatively prime, then the position-sums of the
0-linearizations of a $(k,l)$-arrangement $a$ belong to distinct congruence
classes modulo $l$, where the {\it position-sum} of a 0,1-vector is
the sum of the indices of its 1's.
\PF
We cycle through the 0-linearizations, decreasing the position-sum by $k\mod l$
for each successive 0-linearization.  From one linearization of
$a$, we move to the next by moving the bit in position 1 to position $n=k+l$ and
shifting each other bit down by one.  If we have a 0-linearization with a 0 in
position 1, then a single shift takes us to the next 0-linearization and 
decreases the position-sum by $k$.  If the bit in position 1 is a 1, then
we make additional shifts before moving the first 0 to the back.
For each shift in which a 1 moves from the front to the back, we decrease
the position by one for $k-1$ 1's and increase it by $k+l-1$ for one 1.
The net change in the position-sum is $(k+l-1)-(k-1)=l$.  Thus this operation
does not change the congruence class of the position sum.  Only the last
shift to reach the next 0-linearization changes the congruence class,
again reducing it by $k$ modulo $l$.  \qed

\SH
{5. ({\bigmit k},{\bigmit qk}+{\bigmit p})-ARRANGEMENTS WITH {\bigmit p>}1}
In light of our results about $q$-satisfying sequences of arbitrary lengths,
it is natural to seek comparable extensions of the Strong Cycle Lemma to
$(k,qk+p)$-arrangements.  Unfortunately, when $p>1$ it is possible for
complementary intervals to be $q$-good.  The simple proof of the Strong Cycle
Lemma used when $p=1$ thus fails in the general case, and generalizations for
$p>1$ make weaker statements about $q$-good intervals.  We mention two
special cases of our final theorem: every $0$-linearization of a 
$(k,qk+p)$-arrangement has at least $p$ 0-intervals that are $q$-good, and there
are at least $p$ 0-linearizations in which every $0$-interval is $q$-good.

\TH{9}
(Extended Strong Cycle Lemma)  If $a$ is a $(k,qk+p)$-arrangement and
$p\le i\le qk+p$, then $a$ has at least $qk+2p-i$ 0-linearizations
that have at least $i$ $q$-good 0-intervals.
\PF
The crux of the proof is the {\it augmentation property}: If $b$ is a
0-linearization of a $(k,l)$-arrangement with $l>qk$, and $b'$ is obtained from
$b$ by inserting a 0, then $b'$ has more $q$-good 0-intervals than $b$.  To
prove this, we partition $b$ as $I_1 0 I_2 0 \cdots I_t 0$, where the intervals
$I_j$ contain no ends of $q$-good 0-intervals, and the $t$ elements indicated by
0's are the ends of the $q$-good 0-intervals (some of the $I_j$'s may be empty).
The location of the first $q$-good 0-interval implies that $w_0(I_1)=qw_1(I_1)$.
The location of each successive $q$-good 0-interval implies that
$w_0(I_j)=qw_1(I_j)-1$ for each $j>1$.  To form $b'$, we insert a 0 in some
$I_r$, obtaining $I_r'$.  Each 0 ending a $q$-good 0-interval in $b$ does
so also in $b'$.  In addition, the last 0 in $I_r'$ (which may or may not be
the added 0) also ends a $q$-good 0-interval in $b'$.

We now prove the theorem by induction on $p$.  When $p=1$, the desired statement
is a weakening of the Strong Cycle Lemma.  For $p>1$, we begin by finding a
0-linearization in which every 0-interval is $q$-good.  To do this, delete
$p-1$ 0's arbitrarily to obtain a $(k,qk+1)$-arrangement $\tilde a$.  By the
Strong Cycle Lemma, $\tilde a$ has a unique 0-linearization in which every
0-interval is $q$-good.  Each time we replace one of the deleted 0's, the
augmentation property implies that again every 0-interval is $q$-good.  After
replacing all the deleted 0's, we have a 0-linearization $b$ of $a$ in which
every 0-interval is $q$-qood.

Let $a'$ be the $(k,qk+p-1)$-arrangement obtained by deleting the last element
of $b$ from $a$.  Consider $i$ such that $p\le i\le qk+p$; we have
$p-1 \le i-1 \le qk+p-1$ and $qk+2(p-1)-(i-1)= qk+2p-i-1$.  By the induction
hypothesis, $a'$ has at least $qk+2p-i-1$ 0-linearizations in which at least
$i-1$ 0-intervals are $q$-good.  By the augmentation property, the replacement
of the missing 0 converts these to 0-linearizations of $a$ in which at least
$i$ 0-intervals are $q$-good.  Since every 0-interval in $b$ is $q$-good, $b$
provides the additional needed 0-linearization.  \qed

The extended Strong Cycle Lemma is best possible in the sense that all its
lower bounds may hold with equality simultaneously.  This is achieved by the
$(k,qk+p)$-arrangement in which all the 1's appear together and all the 0's
appear together, which has exactly $p$ 0-linearizations in which all 0-intervals
are $q$-good and one 0-linearization in which exactly $i$ 0-intervals are
$q$-good for each $p\le i<qk+p$.

On the other hand, there may be more 0-linearizations with at least $i$
0-intervals that are $q$-good than guaranteed by the extended Strong Cycle
Lemma, so its inequalities cannot be replaced by equalities.
When $p=tk$, consider the periodic $(k,qk+p)$-arrangement $a$ in which each 1 is
followed by a string of exactly $q+t$ 0's before the next 1.  This arrangement
has exactly $q+t$ ``types'' of 0-linearizations.  When the first 1 in a
0-linearization of $a$ appears after position $q$, every one of the $(q+t)k$
0-intervals is $q$-good; there are $(t+1)k$ such 0-linearizations.  This already
is $k$ more than guaranteed by the Lemma, so the guarantee is exceeded when the
desired number of $q$-good 0-intervals is $i > (q-1)k+p$.  When the first 1
appears in position $j+1$ for some $0\le j \le q$, the number of 0-intervals
that are not $q$-good is $\SE i0r q-j-it$, where $r=\min\{k-1,\FL{(q-j)/t}\}$.
There are $k$ such 0-linearizations for each $j$.  When $k>1$, in this class
of $(k,qk+p)$-arrangements every 0-linearization has more than $p$ 0-intervals
that are $q$-good.

\SH{References}
\frenchspacing
\BP[1]
D.F. Bailey, Counting arrangements of 1's and $-1$'s,
{\it Math. Mag.} 69 (1996), 128--131.
\BP[2]
K.L. Chung and W. Feller, Fluctuations in coin tossing,
{\it Proc. Nat. Acad. Sci. USA} 35 (1949), 605--608.
\BP[3]
N.~Dershowitz and S.~Zaks, The cycle lemma and some applications, \EJC\
11 (1990), 35--40.
\BP[4]
D.A.~Duffus, H.A.~Kierstead, and H.S.~Snevily, An explicit 1-factorization
in the middle of the Boolean lattice, {\it J. Comb. Th. (A)}
65 (1994), 334--342.
\BP[5]
A. Dvoretzky and T. Motzkin, A problem of arrangements,
{\it Duke Math. J.} 14 (1947), 305--313.
\BP[6]
W. Feller, {\it An Introduction to Probability and Its Applications I},
3rd edition.  Wiley \& Sons, 1968.
\BP[7]
R.L.~Graham, D.E.~Knuth, and O.~Patashnik, {\it Concrete Mathematics}.
Addison-Wesley, 1989.
\BP[8]
P.~Hilton and J.~Pedersen, Catalan numbers, their generalization, and
their uses, {\it Math. Intelligencer} 13 (1991), 64--75.
\BP[9]
H.A. Kierstead and W.T. Trotter, Explicit matchings in the middle levels of
the Boolean lattice, {\it Order} 5 (1988), 163--171.
\BP[10]
B. Mont\'agh, A simple proof and a generalization of an old result of Chung
and Feller, \DM\  87 (1991), 105--108.
\BP[11]
T. Narayana, {\it Lattice Path Combinatorics with Statistical Applications}.
Math. Expositions 23, Univ. of Toronto Press, 1979.
\BP[12]
A.M. Odlyzko and H.S. Wilf, The Editor's Corner: $n$ coins in a
fountain, \AMM\ 95 (1988), 840--843.
\BP[13]
G.N. Raney, Functional composition and power series reversion, \TAMS\ 94 (1960),
441--451.
\BP[14]
A.D. Sands, On generalized Catalan numbers, \DM\ 21 (1978), 219--221.
\BP[15]
H.S.~Snevily, {\it Combinatorics of finite sets}.  Ph.D. Thesis,
University of Illinois, 1991.
\bye